\documentclass[12pt]{article}

\usepackage{amsmath,amssymb}

\newtheorem{Remark}{Remark}[section]

\numberwithin{equation}{section}
\newtheorem{Theorem}{Theorem}[section]

\newcommand{\h}{\hspace*{.24in}}

\allowdisplaybreaks
\begin{document}
\title{Overlapping Optimized Schwarz Methods for Parabolic Equations in n-Dimensions}
\author{Minh-Binh TRAN\\
Laboratoire Analyse G\'eom\'etrie et Applications\\
Institut Galil\'ee, Universit\'e Paris 13, France\\
Email: binh@math.univ-paris13.fr}\maketitle
\begin{abstract} 
We introduce in this paper a new tool to prove the convergence of the Overlapping Optimized Schwarz Methods with multisubdomains. The technique is based on some estimates of the errors on the boundaries of the overlapping strips. Our guiding example is an $n$-Dimensional Linear Parabolic Equation.
\end{abstract}
\section{Introduction}
\h In the pioneer work \cite{Lions:1987:OSA}, \cite{Lions:1989:OSA}, \cite{Lions:1990:OSA}, P. L. Lions laid the foundations of the modern theory of Schwarz Methods. With the development of parallel computers, the interest in Schwarz Methods have grown rapidly, as these methods lead to inherently parallel algorithms. However, with Classical Schwarz Methods, high frequency components converge very fast, while low frequency components converge slowly and that slows down the performance of the methods. By replacing Dirichlet Transmission Condition in Classical Schwarz Methods by Robin or higher order Transmission Conditions, we can correct this weekness of Schwarz Method. The new methods are called Optimized Schwarz Methods and have been introduced in \cite{GanderHalpernNataf:1998:OCO}, \cite{GanderHalpernNataf:1999:OSM}. Since then, the convergence properties of the Optimized Schwarz Methods have been studied deeply, based on the following two main tools: Energy Estimates and Laplace and Fourier Transforms. Energy Estimates allow us to study the convergence of the methods in the case of nonoverlapping subdomains. With Energy Estimates, both linear and nonlinear problems have been studied and Optimized Schwarz Methods have been proven to converge, while applying to these equations (see for example, the papers \cite{BenamouDespres:1997:DDM}, \cite{halpern:2009:DGN}, \cite{HalpernSzeftel:2009:NNS}). On the other hand, Laplace and Fourier Transforms allow us to study the convergence of the Overlapping Optimized Schwarz Methods, but for only a few simple equations (see, for example \cite{Bennequin:09:AHB}, \cite{Gander:2007:OSW}, \cite{GanderHalpernNataf:1998:OCO}, \cite{GanderHalpernNataf:1999:OSM}), and the convergence problem of the Domain Decomposition Methods with Robin Transmission Conditions still remains an open problem up to now.
\\\h In this paper, we introduce a new tool to prove the convergence of the Optimized Schwarz Methods for multisubdomains and apply it into an $n$-dimensional linear parabolic equation of the following form
\begin{equation}
\label{1e1}
\frac{\partial u}{\partial t}-\sum_{i,j=1}^na_{i,j}(t)\frac{\partial^2 u}{\partial x_i\partial x_j}+\sum_{i=1}^nb_{i}(t)\frac{\partial u}{\partial x_i}+c(t)u=f(t,x).
\end{equation}
The idea of the technique is to estimate carefully the difference between the values of the errors at the boundaries of the overlapping strips. The technique has the potential to be applied to many other kinds of Partial Differential Equations including nonlinear ones.
\section{Problem Description and Main Results}
We consider the following parabolic equation
\begin{equation}
\label{2e1}
\left \{ \begin{array}{ll}\frac{\partial u}{\partial t}-\sum_{i,j=1}^na_{i,j}(t)\frac{\partial^2 u}{\partial x_i\partial x_j}+\sum_{i=1}^nb_{i}(t)\frac{\partial u}{\partial x_i}+c(t)u=f(t,x),\mbox{ in } (0,T)\times\Omega,\vspace{.1in}\\ 
u(x,t)=g(x,t),\mbox{ on } \partial\Omega\times(0,T),\vspace{.1in}\\ 
u(x,0)=g(x,0),\mbox{ on } \Omega,\end{array}\right. 
\end{equation}
where $\Omega=D\times(\alpha,\beta)$, $D$ is a bounded and smooth enough domain in  $\mathbb{R}^{n-1}$. We impose the following conditions on the coefficients of $(\ref{2e1})$
\\ (A1) For all $i,j$ in $\{1,\dots,I\}$, $a_{i,j}(t)=a_{j,i}(t)$. There exists $\nu_0>0$ such that $A(t)=(a_{i,j}(t))\geq\nu_0I$ for all $t$ belongs to $(0,T)$ in the sense of symmetric positive definite matrices.
\\ (A2) The functions $a_{i,j}$, $b_i$, $c$ are bounded in $C^{\infty}(\mathbb{R})$; $f$ and $g$ are bounded functions in $C^{\infty}(\overline{\Omega\times(0,T)})$.
\\ With the conditions $(A1)$ and $(A2)$, Equation $(\ref{2e1})$ has a unique bounded solution $u$ in $C^{\infty}((0,T)\times\Omega)$. The proof of this result can be infered from Theorems $9$ and $10$, page 71 \cite{Friedman:1964:PDE}. 
\\\h We now divide the domain $\Omega$ into $I$ subdomains, with $\Omega_i=D\times(a_i,b_i)$ and $\alpha=a_1<a_2<b_1<\dots<a_{I}<b_{I-1}<b_I=\beta$. The Optimized Schwarz Waveform Relaxation Algorithm solves $I$ equations in $I$ subdomains instead of solving directly the main problem $(\ref{2e1})$. The iterate $\#k$ in the $l$-th domain, denoted by $u_{l}^{k}$, is defined by
\begin{equation}\label{2e2}
\left \{
\begin{array}{ll}
\frac{\partial u^k_l}{\partial t}-\sum_{i,j=1}^na_{i,j}(t)\frac{\partial^2  u^k_l}{\partial x_i\partial x_j}+\sum_{i=1}^nb_{i}(t)\frac{\partial  u^k_l}{\partial x_i}+c(t) u^k_l=f(t,x)&\mbox{ in }\Omega_l\times(0,T),\vspace{.1in}\\
 \frac{\partial u^k_l(\cdot,a_l,\cdot)}{\partial x_n}+p u^k_l(\cdot,a_l,\cdot)= \frac{\partial u^{k-1}_{l-1}(\cdot,a_l,\cdot)}{\partial x_n}+p u^{k-1}_{l-1}(\cdot,a_l,\cdot) &\mbox{ in }D\times(0,T),\vspace{.1in}\\
 \frac{\partial u^k_l(\cdot,b_l,\cdot)}{\partial x_n}+p u^k_l(\cdot,b_l,\cdot)= \frac{\partial u^{k-1}_{l+1}(\cdot,b_l,\cdot)}{\partial x_n}+p u^{k-1}_{l+1}(\cdot,b_l,\cdot) &\mbox{ in }D\times(0,T),\end{array}\right. 
\end{equation}
here, $p$ is a constant and  for each vector $x$ in $\mathbb{R}^n$, we denote $x=(X,x_n)$, with $X\in\mathbb{R}^{n-1}$ and $x_n\in\mathbb{R}$. 
Each iterate inherits the boundary conditions and the initial values of $u$:
\[
u_l^k(x,t)= g(x,t) \mbox{ on }(\partial\Omega_j\cap\partial\Omega)\times(0,T),\quad
u_l^k(x,0)= g(x,0) \mbox{ in } \Omega_j,
\]
and a special treatment for the extreme subdomains,
\[
u_1^k(\cdot,\alpha,\cdot)=g(\cdot,\alpha,\cdot),
\quad
u_I^k(\cdot,\beta,\cdot)=g(\cdot,\beta,\cdot).
\]
A bounded initial guess $h^0$ in $C^{\infty}(\overline{\Omega\times(0,T)})$ is provided, \textit{i.e.} we solve at step $0$ Equations $(\ref{2e2})$, with boundary data on left and right
\begin{eqnarray*}
\frac{\partial u^1_l(\cdot,a_l,\cdot)}{\partial x_n}+p u^1_l(\cdot,a_l,\cdot)=h^{0}(\cdot,a_l,\cdot) \mbox{ in }D\times(0,T),\nonumber\\
\frac{\partial u^1_l(\cdot,b_l,\cdot)}{\partial x_n}+p u^1_l(\cdot,b_l,\cdot)=h^{0}(\cdot,b_l,\cdot) \mbox{ in }D\times(0,T).
\end{eqnarray*}
By using an induction argument and the same arguments as in Theorem 2, page 144 \cite{Friedman:1964:PDE}, we can see that each subproblem $(\ref{2e2})$ in each iteration has a unique solution. Theorem $10$, page 71 \cite{Friedman:1964:PDE} shows that these solutions belong to $C^{\infty}(\Omega\times(0,T))$. This means that the algorithm is well-posed. 
\\\h Denote by $e_l^k$ the difference between $u_l^k$ and $u$, and substract Equation $(\ref{2e2})$ with the main equation $(\ref{2e1})$, we get the following equation on $e_l^k$
\begin{equation}\label{2e3}
\left \{
\begin{array}{ll}
\frac{\partial e^k_l}{\partial t}-\sum_{i,j=1}^na_{i,j}(t)\frac{\partial^2  e^k_l}{\partial x_i\partial x_j}+\sum_{i=1}^nb_{i}(t)\frac{\partial  e^k_l}{\partial x_i}+c(t) e^k_l=0&\mbox{ in }\Omega_l\times(0,T),\vspace{.1in}\\
 \frac{\partial e^k_l(\cdot,a_l,\cdot)}{\partial x_n}+p e^k_l(\cdot,a_l,\cdot)= \frac{\partial e^{k-1}_{l-1}(\cdot,a_l,\cdot)}{\partial x_n}+p e^{k-1}_{l-1}(\cdot,a_l,\cdot) &\mbox{ in }D\times(0,T),\vspace{.1in}\\
 \frac{\partial e^k_l(\cdot,b_l,\cdot)}{\partial x_n}+p e^k_l(\cdot,b_l,\cdot)= \frac{\partial e^{k-1}_{l+1}(\cdot,b_l,\cdot)}{\partial x_n}+p e^{k-1}_{l+1}(\cdot,b_l,\cdot) &\mbox{ in }D\times(0,T).\end{array}\right. 
\end{equation}
Similarly, each iterate inherits the boundary conditions and the initial values of $u$
\[
e_l^k= 0 \mbox{ on }(\partial\Omega_l\cap\partial\Omega)\times(0,T),\quad
e_l^k(\cdot,\cdot,0)= 0 \mbox{ in } \Omega_l,
\]
and the special treatment for the extreme subdomains,
\[
u_1^k(\cdot,\alpha,\cdot)=0,
\quad
u_I^k(\cdot,\beta,\cdot)=0.
\]
The following theorem states that the algorithm converges.
\begin{Theorem}\label{2t2}Let $\varphi$ be a strictly positive function in $C^1(\mathbb{R})$ such that $$-\max_{x_n\in\mathbb{R}}\left(\frac{\varphi'}{\varphi}(x_n)\right)$$ is large enough, the Optimized Schwarz Waveform Relaxation Method converges in the following sense
$$\lim_{k\to\infty}\max_{l\in\{1,\dots,I\}}\left\Vert\left(\frac{\partial (u_l^k-u)}{\partial x_n}\exp(px_n)\right)^2\varphi(t)\right\Vert_{C(\overline{\Omega_l\times(0,T)})}=0.$$
Moreover, for $l$ in $\{1,\dots,I\}$, the sequence $\{u_l^k-u\}$ converges pointwisely to $0$ as $k$ tends to infinity.
\end{Theorem}
\begin{Remark} We can see that if we choose $\varphi(x_n)=\exp(-\gamma x_n)$, then if $\gamma $ is large enough, $-\max_{x_n\in\mathbb{R}}$ $\left(\frac{\varphi'}{\varphi}(x_n)\right)$ is large enough. The condition of our theorem is then satisfied. 
\end{Remark}
\begin{Remark} Since $a_{i,j}$, $b_i$ are functions of $t$, and the domain is divided in to $n$-subdomains, we cannot use Fourier and Laplace Transforms. Moreover, since the subdomains are overlapping, the Energy Estimates Method cannot be used in our case. In the next section, we introduce a new technique to prove the convergence of the algorithm, the technique is based on the observation that we can estimate the difference between the values of $e^k_l$ on the boundary and in the interior.  
\end{Remark}
\begin{Remark} The result in the theorem remains true if we let $a_{i,j}$, $b_i$ be bounded and continuous functions of $t$ and $x$, but not depend on the $n$-th space variable $x_n$, as we can see in the proof in the following section.
\end{Remark}
\section{The Convergence of the Algorithm}
This section is devoted to the proof of Theorem $\ref{2t2}$. We divide the proof into two steps.
\\{\bf Step 1:} The Error Estimates.
\\\h For $k\in\mathbb{N}$ and $i\in\{1,\dots,I\}$, setting $\epsilon_l^k=e_l^k\exp(px_n)$, we get $e_l^k=\epsilon_l^k\exp(-px_n)$. Equation $(\ref{2e3})$ then leads to
%

\begin{equation}\label{3e1}
\left \{
\begin{array}{ll}
\frac{\partial\epsilon_l^k}{\partial t}-\sum_{i,j=1}^na_{i,j}\frac{\partial^2 \epsilon^k_l}{\partial x_i\partial x_j}+\sum_{i=1}^{n-1}(pa_{i,n}+b_{i})\frac{\partial  \epsilon^k_l}{\partial x_i}+(2pa_{n,n}+b_n)\frac{\partial  \epsilon^k_l}{\partial x_n}\vspace{.1in}\\
~~~~~~~~~~~~~~~~~~~~~~~~~~~~~~~~~~~~~~~~~~+(c-pb_n-p^2a_{n,n})\epsilon^k_l=0,\mbox{ in }\Omega_l\times(0,T), \vspace{.1in}\\
 \frac{\partial \epsilon^k_l(\cdot,a_l,\cdot)}{\partial x_n}= \frac{\partial \epsilon^{k-1}_{l-1}(\cdot,a_l,\cdot)}{\partial x_n} \mbox{ in }D\times(0,T),\vspace{.1in}\\
 \frac{\partial \epsilon^k_l(\cdot,b_l,\cdot)}{\partial x_n}= \frac{\partial \epsilon^{k-1}_{l+1}(\cdot,b_l,\cdot)}{\partial x_n}  \mbox{ in }D\times(0,T),\vspace{.1in}\\
\epsilon_l^k(\cdot,\cdot,\cdot)= 0 \mbox{ on }(\partial\Omega_j\cap\partial\Omega)\times(0,T),\vspace{.1in}\\
\epsilon_l^k(\cdot,\cdot,0)= 0 \mbox{ in } \Omega_j,
\end{array}\right. 
\end{equation}
and for the extreme subdomains,
\[
\epsilon_1^k(\cdot,a_1,\cdot)=0,
\quad
\epsilon_I^k(\cdot,b_I,\cdot)=0.
\]
\h Setting $$\nu_l^k=\frac{\partial\epsilon_l^k}{\partial x_n},$$
we infer from Equation $(\ref{3e1})$ that for $l$ in $\{2,\dots,I-1\}$
\begin{equation}\label{3e2}
\left \{
\begin{array}{ll}
\frac{\partial\nu_l^k}{\partial t}-\sum_{i,j=1}^na_{i,j}\frac{\partial^2 \nu_l^k}{\partial x_i\partial x_j}+\sum_{i=1}^{n-1}(pa_{i,n}+b_{i})\frac{\partial  \nu_l^k}{\partial x_i}+(2pa_{n,n}+b_n)\frac{\partial  \nu_l^k}{\partial x_n}\vspace{.1in}\\
~~~~~~~~~~~~~~~~~~~~~~~~~~~~~~~~~~~~~~~~~~+(c-pb_n-p^2a_{n,n})\nu_l^k=0,\mbox{ in }\Omega_l\times(0,T), \vspace{.1in}\\
 \nu_l^k(\cdot,a_l,\cdot)= \nu^{k-1}_{l-1}(\cdot,a_l,\cdot) \mbox{ in }D\times(0,T),\vspace{.1in}\\
 \nu^k_l(\cdot,b_l,\cdot)= \nu^{k-1}_{l+1}(\cdot,b_l,\cdot) \mbox{ in }D\times(0,T),\vspace{.1in}\\
\nu_l^k(\cdot,\cdot,\cdot)= 0 \mbox{ on }(\partial\Omega_j\cap\partial\Omega)\times(0,T),\vspace{.1in}\\
\nu_l^k(\cdot,\cdot,0)= 0 \mbox{ in } \Omega_j.
\end{array}\right. 
\end{equation}
 We can observe that these systems are with Dicrichlet tramission conditions.
\\\h On $\overline{\Omega_l\times(0,T)}$, we define $$\Phi=\left(\frac{\partial \epsilon_l^k}{\partial x_n}\right)^2\phi(x_n)\varphi(t),$$ where $\phi$ is a strictly positive function in $C^2(\mathbb{R})$ to be chosen later, with the notice that $-\max_{x_n\in\mathbb{R}}\left(\frac{\varphi'}{\varphi}(x_n)\right)$ is large enough. Our purpose is to construct an operator $\mathfrak{L}$ of $\Phi$, such that $\mathfrak{L}(\Phi)$ is negative and then on $\mathfrak{L}$, we can apply the maximum principle to get some estimates on the boundaries for $\Phi$. With these estimates, we can direclty infer some good estimates for $\frac{\partial \epsilon_l^k}{\partial x_n}$ and that lead to our result on the convergence. We now consider the following operator
\begin{eqnarray}\label{3e4}
\mathfrak{H}(\Phi):=\frac{\partial\Phi}{\partial t}-\sum_{i,j=1}^na_{i,j}\frac{\partial^2 \Phi}{\partial x_i\partial x_j}.
\end{eqnarray}
A simple calculation gives
\begin{eqnarray}\label{3e5}
\mathfrak{H}(\Phi)&=&2\nu_l^k\phi\varphi\left(\frac{\partial \nu_l^k}{\partial t}-\sum_{i,j=1}^na_{i,j}\frac{\partial^2\nu_l^k}{\partial x_ix_j}\right)-\sum_{i,j=1}^n2a_{i,j}\phi\varphi\frac{\partial \nu_l^k}{\partial x_i}\frac{\partial \nu_l^k}{\partial x_j}\\\nonumber
& &-\sum_{i=1}^{n}2a_{i,n}\phi'\varphi\nu_l^k\frac{\partial \nu_l^k}{\partial x_i}+\left(\frac{\varphi'}{\varphi}-a_{n,n}\frac{\phi''}{\phi}\right)\phi\varphi(\nu_l^k)^2.
\end{eqnarray}
Since the second term on the right hand side of the previous inequality is negative, it directly leads to 
\begin{eqnarray}\label{3e6}\nonumber
\mathfrak{H}(\Phi)&\leq&2\nu_l^k\phi\varphi\left(\frac{\partial \nu_l^k}{\partial t}-\sum_{i,j=1}^na_{i,j}\frac{\partial^2\nu_l^k}{\partial x_ix_j}\right)\\
& &-\sum_{i=1}^{n}2a_{i,n}\phi'\varphi\nu_l^k\frac{\partial \nu_l^k}{\partial x_i}+\left(\frac{\varphi'}{\varphi}-a_{n,n}\frac{\phi''}{\phi}\right)\phi\varphi(\nu_l^k)^2.\end{eqnarray}
Our purpose is to transfer the right hand side of $(\ref{3e6})$ into the sum of a negative term and a term of $\Phi$, in order to do that, we replace $(\ref{3e2})$ into $(\ref{3e6})$ and get the following bound for $\mathfrak{H}(\Phi)$
\begin{eqnarray}\label{3e7}\nonumber
&&2\nu_l^k\phi\varphi\left(-\sum_{i=1}^{n-1}(pa_{i,n}+b_{i})\frac{\partial  \nu^k_l}{\partial x_i}-(2pa_{n,n}+b_n)\frac{\partial  \nu^k_l}{\partial x_n}-(c-pb_n-p^2a_{n,n})\nu^k_l\right)\\
& &-\sum_{i=1}^{n}2a_{i,n}\phi'\varphi\nu_l^k\frac{\partial \nu_l^k}{\partial x_i}+\left(\frac{\varphi'}{\varphi}-a_{n,n}\frac{\phi''}{\phi}\right)\phi\varphi(\nu_l^k)^2.\end{eqnarray}
Replacing
\begin{eqnarray*}
\frac{\partial\Phi}{\partial x_i}&=&2\phi\varphi\nu_l^k\frac{\partial\nu_l^k}{\partial x_i}, \mbox{ for } i\in\{1,\dots,n-1\},
\end{eqnarray*}
and
\begin{eqnarray*}
\frac{\partial\Phi}{\partial x_n}&=&2\phi\varphi\nu_l^k\frac{\partial\nu_l^k}{\partial x_n}+\phi'\varphi(\nu_l^k)^2,
\end{eqnarray*}
into $(\ref{3e7})$, we get
\begin{eqnarray}\label{3e8}
& &\frac{\partial\Phi}{\partial t}-\sum_{i,j=1}^na_{i,j}\frac{\partial^2 \Phi}{\partial x_i\partial x_j}+\left(2pa_{n,n}+b_n+a_{n,n}\frac{\phi'}{\phi}\right)\frac{\partial  \Phi}{\partial x_n}\\\nonumber
& &+\sum_{i=1}^{n-1}\left(pa_{i,n}+b_{i}+a_{in}\frac{\phi'}{\phi}\right)\frac{\partial  \Phi}{\partial x_i}\\\nonumber
&\leq & \phi\varphi(\nu_l^k)^2\left(\frac{\varphi'}{\varphi}-a_{nn}\frac{\phi''}{\phi}+a_{n,n}\frac{\phi'^2}{\phi^2}-2(c-pb_n-p^2a_{n,n})+(2pa_{n,n}+b_n)\frac{\phi'}{\phi}\right).
\end{eqnarray}
We now get the formula for $\mathfrak{L}$ 
\begin{eqnarray}\label{3e9}
\mathfrak{L}(\Phi)&:=&\frac{\partial\Phi}{\partial t}-\sum_{i,j=1}^na_{i,j}\frac{\partial^2 \Phi}{\partial x_i\partial x_j}+(2pa_{n,n}+b_n+a_{nn}\frac{\phi'}{\phi})\frac{\partial  \Phi}{\partial x_n}\\\nonumber
& &+\sum_{i=1}^{n-1}\left(pa_{i,n}+b_{i}+a_{in}\frac{\phi'}{\phi}\right)\frac{\partial  \Phi}{\partial x_i},
\end{eqnarray}
then if we choose $\varphi$ such that $-\max_{x_n\in\mathbb{R}}\left(\frac{\varphi'}{\varphi}(x_n)\right)$ is large enough, since $a_{n,n}$, $b_n$, $c$, $\frac{\phi'}{\phi}$, $\frac{\phi''}{\phi}$ are all bounded in $C(\mathbb{R})$, we can obtain a negative sign on the right hand side of $(\ref{3e9})$, which means $\mathfrak{L}(\Phi)$ is negative.
\\\h Since $\mathfrak{L}(\Phi)\leq 0$, the maximum of $\Phi$ can only be attained on the boundary of $\Omega_l\times(0,T)$. Since $\Phi=0$ on $\partial\Omega_l\cap\partial\Omega$ and on $\Omega\times\{0\}$, we have the following three estimates.
{\\\h\it Estimate 1}: $1\leq l\leq I$.
\\\h The maximum value(s) of $\Phi$ can be achived on both $\overline{D}\times{\{a_l\}}\times[0,T]$ and $\overline{D}\times{\{b_l\}}\times[0,T]$ and
\begin{eqnarray}\label{3e12}
& &(\nu_l^k(X,x_n,t))^2\phi(x_n)\varphi(t) \leq\\
& \leq &\max\{\max_{\bar D\times[0,T]}\{(\nu_l^k(X,a_{l},t))^2\phi(a_l)\varphi(t)\},\max_{\bar D\times[0,T]}\{(\nu_l^k(X,b_{l},t))^2\phi(b_l)\varphi(t)\}\}\nonumber.
\end{eqnarray}
{\\\h\it Estimate 2}: $l=1$.
\\\h The maximum value(s) of $\Phi$ can be achived on both $\overline{D}\times{\{a_1\}}\times[0,T]$ and $\overline{D}\times{\{b_1\}}\times[0,T]$. If the maximum of $\Phi$ is achived on $\overline{D}\times{\{a_1\}}\times[0,T]$, then at the maximum point, we need that $\frac{\partial\Phi}{\partial n}>0$ due to Hopf's Lemma. We compute
\begin{eqnarray*}
\frac{\partial\Phi}{\partial n}(.,a_1,t)& =&-\frac{\partial\nu_1^k}{\partial x_n}\nu_1^k\phi(a_1)\varphi(t)-(\nu_1^k)^2\phi'(a_1)\varphi(t)\\
& = &-\frac{\partial^2\epsilon_1^k}{\partial x_n^2}\frac{\partial\epsilon_1^k}{\partial x_n }\phi(a_1)\varphi(t)-\left(\frac{\partial\epsilon_1^k}{\partial x_n }\right)^2\phi'(a_1)\varphi(t)\\
& = &-\phi(a_1)\varphi(t)\left[\frac{\partial^2\epsilon_1^k}{\partial x_n^2}\frac{\partial\epsilon_1^k}{\partial x_n }+\left(\frac{\partial\epsilon_1^k}{\partial x_n }\right)^2\frac{\phi'(a_1)}{\phi(a_1)}\right].
\end{eqnarray*}
Since 
\begin{eqnarray*}
\frac{\partial\epsilon_1^k}{\partial t}-\sum_{i,j=1}^na_{i,j}\frac{\partial^2 \epsilon_1^k}{\partial x_i\partial x_j}+\sum_{i=1}^{n-1}(pa_{i,n}+b_{i})\frac{\partial  \epsilon_1^k}{\partial x_i}+(2pa_{n,n}+b_n)\frac{\partial  \epsilon_1^k}{\partial x_n}\\\nonumber
~~~~~~~~~~~~~~~~~~~~~~~~~~~~~~~~~~~~~~~~~~+(c-pb_n-p^2a_{n,n})\epsilon_1^k=0,\mbox{ in }\Omega_1\times(0,T),
\end{eqnarray*}
we get
\begin{eqnarray}\label{3e3}\nonumber
\frac{\partial\epsilon_1^k}{\partial t}(.,a_1,.)-\sum_{i,j=1}^na_{i,j}\frac{\partial^2 \epsilon_1^k}{\partial x_i\partial x_j}(.,a_1,.)+\sum_{i=1}^{n-1}(pa_{i,n}+b_{i})\frac{\partial  \epsilon_1^k}{\partial x_i}(.,a_1,.)\\\nonumber
+(2pa_{n,n}+b_n)\frac{\partial  \epsilon_1^k}{\partial x_n}(.,a_1,.)+(c-pb_n-p^2a_{n,n})\epsilon_1^k(.,a_1,.)=0,\mbox{ in }D\times(0,T).
\end{eqnarray}
Combining this equation and the fact that $\epsilon_l^k(.,a_1,.)=0$ on 
$D\times(0,T)$, we can deduce  
\begin{eqnarray*}
-\frac{\partial^2 \epsilon_1^k}{\partial x_n^2}(.,a_1,.)+\left(2p+\frac{b_n}{a_{n,n}}\right)\frac{\partial  \epsilon_1^k}{\partial x_n}(.,a_1,.)=0,
\end{eqnarray*}
and as a consequence, we can write $\frac{\partial\Phi}{\partial n}$ in a different way
\begin{eqnarray*}
\frac{\partial\Phi}{\partial n}(.,a_1,.)=-\phi(a_1)\varphi(t)\left(\frac{\partial\epsilon_1^k}{\partial x_n }\right)^2\left[\left(2p+\frac{b_n}{a_{n,n}}\right)+\frac{\phi'(a_1)}{\phi(a_1)}\right].
\end{eqnarray*}
With the functions $\phi$ satisfying $$\left(2p+\frac{b_n(t)}{a_{n,n}(t)}\right)+\frac{\phi'(a_1)}{\phi(a_1)}>0,$$ we can see that $$\frac{\partial\Phi}{\partial n}(.,a_1,.)<0;$$ which means that the maximum of $\Phi$ can be achived only on $\overline{D}\times{\{b_1\}}\times[0,T]$, then  
\begin{eqnarray}\label{3e10}
(\nu_{1}^{k}(X,x_n,t))^2\phi(x_n)\varphi(t)\leq\max_{\bar D\times[0,T]}\{(\nu_{1}^{k}(X,b_1,t))^2\phi(b_1)\varphi(t)\}.
\end{eqnarray}
{\\\h\it Estimate 3}: $l=I$.
\\\h The maximum value(s) of $\Phi$ can be achived on both $\overline{D}\times{\{a_I\}}\times[0,T]$ and $\overline{D}\times{\{b_I\}}\times[0,T]$. If the maximum of $\Phi$ is achived on $\overline{D}\times{\{b_I\}}\times[0,T]$, then at the maximum point, we need that $\frac{\partial\Phi}{\partial n}>0$ due to Hopf's Lemma. Similar as in Estimate 2, we can get
\begin{eqnarray*}
\frac{\partial\Phi}{\partial n}(.,b_I,t)& =&\frac{\partial\nu_1^k}{\partial x_n}\nu_1^k\phi(b_I)\varphi(t)+(\nu_1^k)^2\phi'(b_I)\varphi(t)\\
& = &\phi(b_I)\varphi(t)\left[\frac{\partial^2\epsilon_1^k}{\partial x_n^2}\frac{\partial\epsilon_1^k}{\partial x_n }+\left(\frac{\partial\epsilon_1^k}{\partial x_n }\right)^2\frac{\phi'(b_I)}{\phi(b_I)}\right]\\
& = &\phi(b_I)\varphi(t)\left(\frac{\partial\epsilon_1^k}{\partial x_n }\right)^2\left[\left(2p+\frac{b_n}{a_{n,n}}\right)+\frac{\phi'(b_I)}{\phi(b_I)}\right].
\end{eqnarray*}
With the functions $\phi$ satisfying 
$$\left(2p+\frac{b_n(t)}{a_{n,n}(t)}\right)+\frac{\phi'(b_I)}{\phi(b_I)}<0,$$ we can see that $$\frac{\partial\Phi}{\partial n}(.,b_I,.)<0;$$ which means the maximum of $\Phi$ can be achived only on $\overline{D}\times{\{a_I\}}\times[0,T]$, then  
\begin{eqnarray}\label{3e11}
 (\nu_{I}^{k}(X,x_n,t))^2\phi(x_n)\varphi(t)\leq \max_{\bar D\times[0,T]}\{(\nu_{I}^{k}(X,a_{I},t))^2\phi(a_I)\varphi(t)\}.
\end{eqnarray}
{\h\bf Step 2:} Proof of convergence, $$\lim_{k\to\infty}\max_{l\in\left\{1,\dots,I\right\}}\left\Vert\left(\frac{\partial\epsilon_l^k}{\partial x_n}\right)^2\varphi(t)\right\Vert_{C(\overline{\Omega_l\times(0,T)})}=0.$$
\\\h In the proof of convergence, we will use the three estimates $(\ref{3e12})$, $(\ref{3e10})$ and $(\ref{3e11})$ by fixing the function $\varphi(t)$ and replacing $\phi(x_n)$ by appropriate functions $\overline\phi_i$, $\tilde\phi_i$, $\overline\phi_*$, $\tilde\phi_*$ $(i\in\{1,\dots,I\})$ in each subdomain. Before coming to the details of the proof, we need a notation
$$E_k=\max_{l\in\{1,\dots,I\}}||(\nu^k_l)^2\phi\varphi||_{C(\overline{\Omega_l\times(0,T)})}.$$
{\\\h\bf Step 2.1:} Estimate of the right boudaries of the sub-domains.
\\\h Consider the $I$-th domain, at the $k$-th step, $(\ref{3e11})$ infers
$$(\nu_I^k(X,x_n,t))^2\overline\phi_I(x_n)\varphi(t)\leq \max_{\bar D\times[0,T]}\{(\nu_I^k(X,a_{I},t))^2\overline\phi_I(a_{I})\varphi(t)\},$$
where $\overline\phi_I$ is a strictly postive function and will be chosen later.
\\ Replace $x_n$ by $b_{I-1}$, we get
$$(\nu_I^k(X,b_{I-1},t))^2\overline\phi_I(b_{I-1})\varphi(t)\leq \max_{\bar D\times[0,T]}\{(\nu_I^k(X,a_{I},t))^2\overline\phi_I(a_{I})\varphi(t)\}.$$
Since $\nu_I^k(X,b_{I-1},t)=\nu_{I-1}^{k+1}(X,b_{I-1},t)$, 
$$(\nu_{I-1}^{k+1}(X,b_{I-1},t))^2\overline\phi_I(b_{I-1})\varphi(t)\leq \max_{\bar D\times[0,T]}\{(\nu_I^k(X,a_{I},t))^2\overline\phi_I(a_{I})\varphi(t)\}.$$
The inequality becomes
$$(\nu_{I-1}^{k+1}(X,b_{I-1},t))^2\varphi(t)\leq\frac{\overline\phi_I(a_{I})}{\overline\phi_I(b_{I-1})}\max_{\bar D\times[0,T]}\{(\nu_I^k(X,a_{I},t))^2\varphi(t)\}.$$
We can choose $\overline\phi_I$ such that $\frac{\overline\phi_I(a_{I})}{\overline\phi_I(b_{I-1})}<1$, and deduce 
\begin{equation}\label{3e13}
(\nu_{I-1}^{k+1}(X,b_{I-1},t))^2\varphi(t)\leq\frac{\overline\phi_I(a_{I})}{\overline\phi_I(b_{I-1})}E_k.
\end{equation}
\h Moreover, on the $(I-1)$-th domain, at the $(k+1)$-th step, $(\ref{3e12})$ leads to
\begin{eqnarray*}
& &(\nu_{I-1}^{k+1}(X,x_n,t))^2\overline\phi_{I-1}(x_n)\varphi(t) \leq \\
 & &\max\{\max_{\bar D\times[0,T]}\{(\nu_{I-1}^{k+1}(X,b_{I-1},t))^2\overline\phi_{I-1}(b_{I-1})\varphi(t)\},\\
& &\max_{\bar D\times[0,T]}\{(\nu_{I-1}^{k+1}(X,a_{I-1},t))^2\overline\phi_{I-1}(a_{I-1})\varphi(t)\}\},
\end{eqnarray*}
where $\overline\phi_{I-1}$ is a strictly postive function and will be chosen later.
\\ Since $\nu_{I-1}^{k+1}(X,b_{I-2},t)=\nu_{I-2}^{k+2}(X,b_{I-2},t)$, 
\begin{eqnarray*}
& &(\nu_{I-2}^{k+2}(X,b_{I-2},t))^2\overline\phi_{I-1}(b_{I-2})\varphi(t) \leq \\
 & &\max\{\max_{\bar D\times[0,T]}\{(\nu_{I-1}^{k+1}(X,b_{I-1},t))^2\overline\phi_{I-1}(b_{I-1})\varphi(t)\},\\
& &\max_{\bar D\times[0,T]}\{(\nu_{I-1}^{k+1}(X,a_{I-1},t))^2\overline\phi_{I-1}(a_{I-1})\varphi(t)\}\}.
\end{eqnarray*}
Hence
\begin{eqnarray*}
& &(\nu_{I-2}^{k+2}(X,b_{I-2},t))^2\varphi(t) \leq \\
 & &\max\left\{\frac{\overline\phi_{I-1}(b_{I-1})}{\overline\phi_{I-1}(b_{I-2})}\max_{\bar D\times[0,T]}\{(\nu_{I-1}^{k+1}(X,b_{I-1},t))^2\varphi(t)\right.\},\\
& &\max_{\bar D\times[0,T]}\left.\left\{\frac{\overline\phi_{I-1}(a_{I-1})}{\overline\phi_{I-1}(b_{I-2})}(\nu_{I-1}^{k+1}(X,a_{I-1},t))^2\varphi(t)\right\}\right\}.
\end{eqnarray*}
Combining this inequality with $(\ref{3e13})$, we get
\begin{eqnarray*}
(\nu_{I-2}^{k+2}(X,b_{I-2},t))^2\varphi(t) \leq \max\left\{\frac{\overline\phi_{I-1}(b_{I-1})}{\overline\phi_{I-1}(b_{I-2})}\frac{\overline\phi_I(a_{I})}{\overline\phi_I(b_{I-1})}E_k,\frac{\overline\phi_{I-1}(a_{I-1})}{\overline\phi_{I-1}(b_{I-2})}E_{k+1}\right\}.
\end{eqnarray*}
We choose $\overline\phi_{I-1}$ such that $$\frac{\overline\phi_{I-1}(b_{I-1})}{\overline\phi_{I-1}(b_{I-2})}\frac{\overline\phi_I(a_{I})}{\overline\phi_I(b_{I-1})}=\frac{\overline\phi_{I-1}(a_{I-1})}{\overline\phi_{I-1}(b_{I-2})}<1,$$ which is equivalent to $$\frac{\overline\phi_I(a_{I})}{\overline\phi_I(b_{I-1})}=\frac{\overline\phi_{I-1}(a_{I-1})}{\overline\phi_{I-1}(b_{I-2})}<1.$$ As a consequence, 
\begin{equation}\label{3e14}
(\nu_{I-2}^{k+2}(X,b_{I-2},t))^2\varphi(t)\leq\frac{\overline\phi_{I-1}(a_{I-1})}{\overline\phi_{I-1}(b_{I-2})}\max\{E_{k},E_{k+1}\}.
\end{equation}
\h Using the same techniques as the ones we use to achive $(\ref{3e13})$ and $(\ref{3e14})$, we can prove that
\begin{equation}\label{3e15}
(\nu_{I-j}^{k+j}(X,b_{I-j},t))^2\varphi(t)\leq\frac{\overline\phi_{I-j+1}(a_{I-j+1})}{\overline\phi_{I-j+1}(b_{I-j})}\max\{E_{k},\dots,E_{k+j-1}\}, 
\end{equation}
where $\overline\phi_{I-j+1}$ is a strictly positive function satisfying 
$$\frac{\overline\phi_{I-j+1}(a_{I-j+1})}{\overline\phi_{I-j+1}(b_{I-j})}<1,$$
with $j=\{1,\dots,I-1\}$.
\\\h Now, with $(\ref{3e10})$, we can choose a strictly postive function $\overline\phi_{*}$ such that $\overline\phi_{*}(b_I)>\overline\phi_*(a_I)$, then
$$(\nu_I^k(X,b_I,t))^2\overline\phi_*(b_I)\varphi(t)\leq \max_{\bar D\times[0,T]}\{(\nu_I^k(X,a_{I},t))^2\overline\phi_*(a_{I})\varphi(t)\},$$
and as a result
$$(\nu_I^k(X,b_I,t))^2\overline\phi_*(b_I)\varphi(t)\leq \max_{\bar D\times[0,T]}\{(\nu_{I-1}^{k-1}(X,a_{I},t))^2\overline\phi_*(a_{I})\varphi(t)\},$$
which implies
\begin{equation}\label{3e16}
(\nu_I^k(X,b_I,t))^2\varphi(t)\leq \frac{\overline\phi_*(a_{I})}{\overline\phi_*(b_I)}E_{k-1}.
\end{equation}
{\h\bf Step 2.2:} Estimate of the left boundaries of the sub-domains.
\\\h Consider the $1$-th domain, at the $k$-th step, $(\ref{3e11})$ infers
$$(\nu_1^k(X,x_n,t))^2\tilde\phi_1(x_n)\varphi(t)\leq \max_{\bar D\times[0,T]}\{(\nu_1^k(X,b_{1},t))^2\tilde\phi_1(b_{1})\varphi(t)\},$$
where $\tilde\phi_1$ is a strictly postive function and will be chosen later.\\
Replacing $x_n$ by $a_{2}$, we get
$$(\nu_1^k(X,a_{2},t))^2\varphi(t)\leq \frac{\tilde\phi_1(b_{1})}{\tilde\phi_1(a_{2})}\max_{\bar D\times[0,T]}\{(\nu_1^k(X,b_{1},t))^2\varphi(t)\}.$$
Since $\nu_1^k(X,a_{2},t)=\nu_2^{k+1}(X,a_{2},t)$, then
$$(\nu_2^{k+1}(X,a_{2},t))^2\varphi(t)\leq \frac{\tilde\phi_1(b_{1})}{\tilde\phi_1(a_{2})}\max_{\bar D\times[0,T]}\{(\nu_1^k(X,b_{1},t))^2\varphi(t)\}.$$
We choose $\tilde\phi_1$ such that $$\frac{\tilde\phi_1(b_{1})}{\tilde\phi_1(a_{2})}<1,$$ and deduce
\begin{equation}\label{3e17}
(\nu_2^{k+1}(X,a_{2},t))^2\varphi(t)\leq \frac{\tilde\phi_1(b_{1})}{\tilde\phi_1(a_{2})}E_k.
\end{equation}
\h Moreover, on the $2$-th domain, at the $(k+1)$-th step, $(\ref{3e11})$ leads to
\begin{eqnarray*}
& &(\nu_{2}^{k+1}(X,x_n,t))^2\tilde\phi_2(x_n)\varphi(t)\leq\\
& & \max\{\max_{\bar D\times[0,T]}\{(\nu_{2}^{k+1}(X,b_{2},t))^2\tilde\phi_2(b_2)\varphi(t)\},\\
& &\max_{\bar D\times[0,T]}\{(\nu_{2}^{k+1}(X,a_{2},t))^2\tilde\phi_2(a_2)\varphi(t)\}\},
\end{eqnarray*}
where $\tilde\phi_2$ is a strictly postive function and will be chosen later.\\
Since $\nu_{2}^{k+1}(X,a_{3},t)=\nu_{3}^{k+2}(X,a_{3},t)$, then
\begin{eqnarray*}
& &(\nu_{3}^{k+2}(X,a_3,t))^2\tilde\phi_2(a_3)\varphi(t)\leq\\
& & \max\{\max_{\bar D\times[0,T]}\{(\nu_{2}^{k+1}(X,b_{2},t))^2\tilde\phi_2(b_2)\varphi(t)\},\\
& &\max_{\bar D\times[0,T]}\{(\nu_{2}^{k+1}(X,a_{2},t))^2\tilde\phi_2(a_2)\varphi(t)\}\}.
\end{eqnarray*}
Hence
\begin{eqnarray*}
& &(\nu_{3}^{k+2}(X,a_3,t))^2\tilde\varphi(t)\leq\\
& & \max\left\{\frac{\tilde\phi_2(b_2)}{\tilde\phi_2(a_3)}\max_{\bar D\times[0,T]}\{(\nu_{2}^{k+1}(X,b_{2},t))^2\varphi(t)\}\right.,\\
& &\left.\frac{\tilde\phi_2(a_2)}{\tilde\phi_2(a_3)}\max_{\bar D\times[0,T]}\{(\nu_{2}^{k+1}(X,a_{2},t))^2\varphi(t)\}\right\}.
\end{eqnarray*}
Combining this with $(\ref{3e17})$, we get
\begin{eqnarray*}
(\nu_{3}^{k+2}(X,a_3,t))^2\varphi(t)\leq
 \max\left\{\frac{\tilde\phi_2(b_2)}{\tilde\phi_2(a_3)}E_{k+1},\frac{\tilde\phi_2(a_2)}{\tilde\phi_2(a_3)}\frac{\tilde\phi_1(b_{1})}{\tilde\phi_1(a_{2})}E_k\right\}.
\end{eqnarray*}We choose $\tilde\phi_2$ such that $$\frac{\tilde\phi_2(b_2)}{\tilde\phi_2(a_3)}=\frac{\tilde\phi_2(a_2)}{\tilde\phi_2(a_3)}\frac{\tilde\phi_1(b_{1})}{\tilde\phi_1(a_{2})}<1,$$ which is equivalent to $$\frac{\tilde\phi_2(b_2)}{\tilde\phi_2(a_2)}=\frac{\tilde\phi_1(b_{1})}{\tilde\phi_1(a_{2})}<1.$$ We then obtain
\begin{equation}\label{3e18}
(\nu_{3}^{k+2}(X,a_3,t))^2\varphi(t)\leq\frac{\tilde\phi_2(b_2)}{\tilde\phi_2(a_3)}\max\{E_{k},E_{k+1}\}.
\end{equation}
\h Using the same techniques as the ones that we use to achive $(\ref{3e17})$ and $(\ref{3e18})$, we can prove that
\begin{equation}\label{3e19}
(\nu_{j}^{k+j-1}(X,a_{j},t))^2\varphi(t)\leq\frac{\tilde\phi_{j-1}(b_{j-1})}{\tilde\phi_{j-1}(a_{j})}\max\{E_{k},\dots,E_{k+j-2}\}, 
\end{equation}
where $\tilde\phi_{j-1}$ is a stricly positive function satisfying
$$\frac{\tilde\phi_{j-1}(b_{j-1})}{\tilde\phi_{j-1}(a_{j})}<1,$$
with $j$ belongs to $\{1,\dots,I-1\}$.
\\\h Now, with $(\ref{3e11})$, we can choose a strictly positive function $\tilde\phi_{*}$ such that $\tilde\phi_{*}(b_1)<\tilde\phi_*(a_1)$, and get
$$(\nu_1^k(X,a_1,t))^2\tilde\phi_*(a_1)\varphi(t)\leq \max_{\bar D\times[0,T]}\{(\nu_1^k(X,b_{1},t))^2\tilde\phi_*(b_{1})\varphi(t)\},$$
which is equivalent to
$$(\nu_1^k(X,a_1,t))^2\tilde\phi_*(a_1)\varphi(t)\leq \max_{\bar D\times[0,T]}\{(\nu_2^{k-1}(X,b_{1},t))^2\tilde\phi_*(b_{1})\varphi(t)\}.$$
That implies
\begin{equation}\label{3e20}
(\nu_I^k(X,a_1,t))^2\varphi(t)\leq \frac{\tilde\phi_*(b_{1})}{\tilde\phi_*(a_1)}E_{k-1}.
\end{equation}
{\h\bf Step 2.3:} Convergence result.
\\\h From $(\ref{3e15})$, $(\ref{3e16})$, $(\ref{3e19})$ and $(\ref{3e20})$, there exists $\gamma$ in $(0,1)$ such that
\begin{eqnarray}\label{3e21}
(\nu_l^{k+I}(X,a_l,t))^2\varphi(t)\leq\gamma\max\{E_{k},\dots,E_{k+I-1}\}, \mbox{ for }l\in\{1,\dots,I\},
\end{eqnarray}
and
\begin{eqnarray}\label{3e22}
(\nu_l^{k+I}(X,b_l,t))^2\varphi(t)\leq\gamma\max\{E_{k},\dots,E_{k+I-1}\}, \mbox{ for }l\in\{1,\dots,I\}. 
\end{eqnarray}
\h Using $(\ref{3e12})$ for $\phi\equiv1$, we have that 
\begin{eqnarray}\label{3e23}
& &(\nu_l^k(X,x_n,t))^2\varphi(t) \leq\\
& \leq &\max\{\max_{\bar D\times[0,T]}\{(\nu_l^k(X,a_{l},t))^2\varphi(t)\},\max_{\bar D\times[0,T]}\{(\nu_l^k(X,b_{l},t))^2\varphi(t)\}\}\nonumber.
\end{eqnarray}
\h Combining $(\ref{3e21})$, $(\ref{3e22})$ and $(\ref{3e23})$, we get
\begin{eqnarray}\label{3e24}
E_{k+I}\leq\gamma\max\{E_{k},\dots,E_{k+I-1}\}.
\end{eqnarray}
\h Hence, $E_k$ tends to $0$ as $k$ tends to infinity. That gives 
$$\lim_{k\to\infty}\max_{l\in\{1,\dots,I\}}\left\Vert\left(\frac{\partial\epsilon_l^k}{\partial x_n}\right)^2\varphi(t)\right\Vert_{C(\overline{\Omega_l\times(0,T)})}=0.$$
{\h\bf Step 3:} Proof of convergence: for $l$ in $\{1,\dots,I\}$, the sequence $\{e_l^k\}$ converges pointwisely to $0$ as $k$ tends to infinity.
\\\h Since for $l$ in $\{1,\dots,I\}$, $$\lim_{k\to\infty}\left\Vert\left(\frac{\partial\epsilon_l^k}{\partial x_n}\right)^2\varphi(t)\right\Vert_{C(\overline{\Omega_l\times(0,T)})}=0,$$ then $\left\{\left(\frac{\partial\epsilon_l^k}{\partial x_n}\right)^2\varphi(t)\right\}$ converges to $0$ pointwisely and the sequence is bounded by a constant $M_0$. Since $\varphi$ is strictly positive on $[0,T]$, there exist positive constants $M_1,$ $M_2$ such that $M_1<\varphi<M_2$. That means the sequence $\left\{\left|\frac{\partial\epsilon_l^k}{\partial x_n}\right|\right\}$ converges to $0$ pointwisely  and is bounded by a constant $M_3$. 
\\\h For $l=1$, with a fixed value of $(X,t)$, the Lebesgue Dominated Convergence Theorem confirms that for $x_n$ belongs to $[a_1,b_1]$, $$\int_{a_1}^{x_n}\frac{\partial\epsilon_1^k}{\partial x_n}(X,\zeta,t)d\zeta$$ converges to $0$ as $k$ tends to infinity. Hence the sequence $$\{\epsilon_1^k(X,x_n,t)-\epsilon_1^k(X,a_1,t)\}$$ converges to $0$ as $k$ tends to infinity. Since $\epsilon_1^k(X,a_1,t)=0$, then the sequence $\{\epsilon_1^k\}$ converges to $0$ pointwisely.
\\\h For $l=2$, with a fixed value of $(X,t)$, again by the Lebesgue Dominated Convergence Theorem, for $x_n$ belongs to $[a_2,b_2]$, the sequence $$\left\{\int_{a_2}^{x_nt}\frac{\partial\epsilon_2^k}{\partial x_n}(X,\zeta,t)d\zeta\right\}$$ converges to $0$ as $k$ tends to infinity. Hence the sequence $$\left\{\epsilon_2^k(X,x_n,t)-\epsilon_2^k(X,a_2,t)\right\}$$ converges to $0$ as $k$ tends to infinity. Since $$\frac{\partial e_2^k}{\partial x_n}(X,a_2,t)+pe_2^k(X,a_2,t)=\frac{\partial e_1^{k-1}}{\partial x_n}(X,a_2,t)+pe_1^{k-1}(X,a_2,t),$$ and the sequences $\{e_1^k\}$, $\left\{\left|\frac{\partial\epsilon_l^k}{\partial x_n}\right|\right\}$ converge to $0$ pointwisely for $l$ belongs to $\{1,\dots,I\}$, we can deduce that $\epsilon_2^k(X,x_n,t)$ converges to $0$ as $k$ tends to infinity. 
\\\h By similar processes, we can prove that for $l$ in $\{1,\dots,I\}$, the sequence $\{e_l^k\}$ converges pointwisely to $0$ as $k$ tends to infinity. This concludes the proof.
{\\\bf Acknowledgement.} The author would like to thank Professor Laurence Halpern for her kindness and support.
\bibliographystyle{plain}\bibliography{ODDMforParabolic}
\end{document}